\documentclass[11pt]{article}

\usepackage[utf8]{inputenc}
\usepackage{graphicx}
\usepackage{amssymb}
\usepackage{amsmath}
\usepackage{amsthm}
\usepackage{color}
\usepackage{tikz}
\usepackage{todonotes}
\usepackage{subcaption}
\usepackage{float}
\usepackage[justification=centering]{caption}
\usepackage[
    natbib=true,
    style=alphabetic,
    sorting=nyt
]{biblatex}
\newtheorem{thm}{Theorem}[section]

\theoremstyle{definition}
\newtheorem{defn}[thm]{Definition}
\newtheorem{example}[thm]{Example}
\newtheorem{conjecture}[thm]{Conjecture}
\newtheorem{remark}[thm]{Remark}
\newtheorem{question}[thm]{Question}
\newtheorem{notation}[thm]{Notation}

\usepackage{hyperref}
\usepackage{comment}

\newcommand\E{\mathbb{E}}

\newcommand{\F}{\mathbb{F}}

\renewcommand\P{\mathcal{P}}

\newcommand\Z{\mathbb{Z}}
\newcommand{\C}{\mathbb{C}}
\newcommand{\bbC}{\mathbb{C}}

\renewcommand{\H}{\mathbb{H}}
\newcommand{\A}{\mathcal{A}}

\newcommand{\bnd}{\partial}
\newcommand{\inverse}{^{-1}}

\DeclareMathOperator{\cat}{CAT}
\DeclareMathOperator{\PSL}{PSL}

\renewcommand\hat{\widehat}

\title{
	Hyperbolic boundaries vs. hyperbolic groups}
\author{Michael Ben-Zvi, Jiayi Lou, and Genevieve S. Walsh}

\date{}
\begin{document}

	\maketitle

	\abstract{ The aim of these notes is to connect the theory of hyperbolic and relatively hyperbolic groups to the theory of manifolds and Kleinian groups. We also give definitions and many examples of relatively hyperbolic groups and their boundaries. We survey some of the extensive work that has been done in the field. These notes are based on lectures given by the third author at CIRM in the Summer of 2018.}

\section{Introduction and Preliminaries}
 From the three-manifold theorist's point of view, hyperbolic and relatively hyperbolic groups are generalizations of Kleinian groups. Here we highlight some deep connections between the two theories. All groups are assumed to be finitely generated and all manifolds irreducible and orientable, unless otherwise specified.

 Trees, $\mathbb{H}^2$, and $\mathbb{H}^n$ are all examples of hyperbolic metric spaces.  Similarly, free groups, the fundamental groups of closed hyperbolic surfaces and the fundamental groups of closed hyperbolic three-manifolds, are examples of hyperbolic groups. The fundamental groups of cusped hyperbolic 3-manifolds, such as hyperbolic knot groups, are not hyperbolic groups.  They are however, relatively hyperbolic groups, as are all geometrically finite Kleinian groups. We recall the definition of hyperbolic metric spaces and groups (the reader should verify the examples above).

\begin{defn}

A \textit{hyperbolic metric} space is a geodesic metric space such that geodesic triangles are slim.  That is, there is a global constant $\delta > 0$ such that for all geodesic triangles the third side is contained in the $\delta$ neighborhood of the union of the other two.  A group acts {\it geometrically} on a proper metric space if the action is properly discontinuous, isometric and co-compact.
A {\it hyperbolic group} is a group which acts geometrically on {\it some} proper hyperbolic metric space $X$.
\end{defn}

 A canonical example is a {\it co-compact Kleinian group}, a group which acts geometrically on $\mathbb{H}^3$.

Hyperbolic groups are often called {\it Gromov} hyperbolic groups after \cite{Gromov87}. We also direct the reader to the several excellent surveys: \cite{ABC91}, \cite{BowditchGGTnotes}, \cite{CDP90}, among others.

\begin{notation} Throughout these notes we will discuss several boundaries for metric spaces, and use notation as follows.
\begin{itemize}
\item $\partial X$: We will denote the Gromov boundary of a hyperbolic space by $\partial X$.  Similarly, we will denote the visual boundary of a $\cat(0)$ space by $\partial X$.  We hope no confusion arises here.  When a geodesic metric space is both $\cat(0)$ and hyperbolic, the boundaries are homeomorphic. For definitions of these boundaries and this fact see \cite{BH}.

\item $\partial G$: We denote the Gromov boundary of a hyperbolic group by $\partial G$, which is the boundary of any hyperbolic space on which $G$ acts geometrically. This is well-defined since  any two spaces $X$ and $Y$ that $G$ acts upon geometrically are quasi-isometric, which implies that $\partial X$ and $\partial Y$ are homeomorphic.  When the boundary of a $\cat(0)$ group is well-defined, we will use the same notation $\partial G$. While it is known that there are many examples of $\cat(0)$ groups that do not have well-defined boundaries \cite{Croke-Kleiner}, we will entirely restrict ourselves to $\cat(0)$ groups with isolated flats (see Definition \ref{IsolatedFlats}), which do have well-defined boundaries \cite{HK05}.

\item $\partial(G, \mathcal{P})$: This will denote the Bowditch boundary of the relatively hyperbolic pair $(G, \mathcal{P})$. For more elaboration, see Definition \ref{RelHypBnd}.
\end{itemize}
\end{notation}

\noindent {\bf Plan of paper:} In Section \ref{day1} we discuss relatively hyperbolic groups and their boundaries, making relations with the hyperbolic and $\cat(0)$ visual boundaries.  We also give several examples.
In Section \ref{day2} we discuss equivalent definitions of relative hyperbolicity, and some spaces that are useful.  In section \ref{day3} we discuss the connection with Kleinian groups, and also some algebraic information that can be gleaned from these boundaries.

\noindent {\bf Acknowledgements.} We thank CIRM for hosting the school where these lectures originated, and the referee for very helpful comments which  improved this paper. The third author was partially supported by NSF grant DMS - 2005353.
\section{Relatively hyperbolic groups and their boundaries}
\label{day1}

 A {\it geometrically finite Kleinian group} is a Kleinian group which acts geometrically finitely on the convex hull of its limit set.  See Section \ref{day3} for more detailed definitions. More generally,  a {\it relatively hyperbolic group pair $(G, \mathcal{P})$} is a group pair that acts geometrically finitely on a proper hyperbolic metric space $X$.

There are many equivalent definitions of geometrically finite Kleinian groups, (see Bowditch \cite{BowGF}). Similarly, there are many equivalent definitions of relatively hyperbolic groups.  We will take as our definition \cite[Def 2]{Bowrelhyp}.  Other, equivalent definitions are discussed in Section \ref{day2}.

First, we define conical limit points and bounded parabolic points.

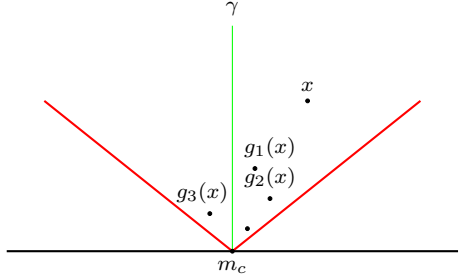
\begin{figure}
    \centering
    \begin{tikzpicture}

    \draw[thick,red] (-2.5,2) -- (0,0);
    \draw[thick,red] (2.5,2) -- (0,0);
    \draw[green] (0,3) -- (0,0);
    \draw[thick] (-3,0) -- (3,0);

    \node[above] at (0,3) {\scriptsize$\gamma$};
    \node[below] at (0,0) {\scriptsize$m_c$};

    \node[above] at (1,2) {\scriptsize$x$};
    \draw[fill] (1,2) circle [radius=0.025];

     \node[above] at (.3,1.1) {\scriptsize$~~~~g_1(x)$};
    \draw[fill] (.3,1.1) circle [radius=0.025];

     \node[above] at (.5,.7) {\scriptsize$g_2(x)$};
    \draw[fill] (.5,.7) circle [radius=0.025];

    \node[above] at (-.4,.5) {\scriptsize $g_3(x)$};
    \draw[fill] (-.3,.5) circle [radius=0.025];

    \draw[fill] (.2,.3) circle [radius=0.025];

    \draw[fill] (0,0) circle [radius=0.025];

    \end{tikzpicture}
    \caption{A conical limit point}
    \label{fig:conical}
\end{figure}

\begin{defn}
		Let $X$ be a proper, hyperbolic geodesic metric space where $G$ acts on $X$ properly discontinuously by isometries.

	\begin{itemize}

		\item Conical limit point: A point $m\in \bnd X$ is a \textit{conical limit point} if there is a geodesic $\gamma\to m$, a point $x\in X$, and a sequence of elements $\{g_n\} \subset G$ such that $g_nx \to m$ and $d(g_nx,\gamma) < r$ for some $r>0$. See Figure \ref{fig:conical}.

		\item Parabolic: $P\leq G$ is a \textit{parabolic subgroup} if it is infinite, contains no loxodromic elements, and fixes a point $x_P\in \bnd X$. In this case, $x_P$ is called a \textit{parabolic point}.

		\item Bounded parabolic: A parabolic point $x_P\in \bnd X$ is \textit{bounded} if $\bnd X \setminus \{x_P\} / P$ is compact.
	\end{itemize}
\end{defn}

Note that if an infinite subgroup fixes more than one point it must contain a loxodromic element, so a parabolic subgroup must fix exactly 1 point.

\begin{defn}[\cite{Bowrelhyp}] \label{def:original rel hyp} A group pair is a group $G$ and a family $\P$ of infinite subgroups consisting of finitely many conjugacy classes.  The pair $(G,\P)$ is \textit{relatively hyperbolic} if $G$ acts on $X$ properly discontinuously and by isometries, where $X$ is a proper hyperbolic geodesic metric space such that:

\begin{enumerate}
	\item each point of $\bnd X$ is either a conical limit point or a bounded parabolic point.
	\item $\P$ is exactly the collection of maximal parabolic subgroups.
\end{enumerate}

In the case that we have a properly discontinuous action by isometries and these two conditions are satisfied, we say $(G, \mathcal{P})$ acts \textit{geometrically finitely} on $X$.

The elements of $\P$ are called \textit{peripheral} subgroups.
\end{defn}

\begin{defn} \label{RelHypBnd} The relatively hyperbolic boundary $\partial(G, \mathcal{P})$ (alternatively  $\partial_B(G, \mathcal{P})$) of the group pair $(G, \P)$ is the boundary of any hyperbolic metric space that $(G, \mathcal{P})$ acts on geometrically finitely. This is also called the Bowditch boundary.
\end{defn}

Two such spaces are not necessarily equivariantly quasi-isometric \cite{Burnsrigid}, as in the hyperbolic case.  However, the relatively hyperbolic boundary is still well-defined up to homeomorphism, \cite[Section 9]{Bowrelhyp}. In the case of a hyperbolic knot complement such as the figure-eight knot complement, the Bowditch boundary of the group pair is $S^2$, the boundary of $\mathbb{H}^3$.  It is an open question to understand relatively hyperbolic group pairs with Bowditch boundary $S^2$ or even a subset of $S^2$ (those with planar boundary).  Variants of this question have been explored in \cite{GMS, TW, HW1,MartinSkora89}. See Question \ref{question: vgfk} and Conjecture \ref{conjecture: vgfk} for further discussion.

There are lots of natural relatively hyperbolic groups.  One good source of examples is the class of hyperbolic groups. Given an almost malnormal collection of quasiconvex subgroups (which might be the empty set), one obtains a relatively hyperbolic group pair.  Let $G$ be a hyperbolic group and $\P$ a collection of quasiconvex subgroups. We say the collection is \textit{almost malnormal} if for every $P,P'\in \P$ and $g\in G$, whenever

        $$|P' \cap gPg\inverse| = \infty$$
        then $P=P'$ and $g\in P$.

\begin{thm}[{\cite[Theorem 7.11]{Bowrelhyp}}]
Let $G$ be a non-elementary hyperbolic group and $\P$ an almost malnormal collection of infinite quasiconvex subgroups of $G$, where $\P$ is invariant under conjugacy and contains finitely many conjugacy classes. Then $(G,\P)$ is relatively hyperbolic.
\end{thm}

\begin{example}\label{example:ap} To illustrate how the Bowditch boundary can change dramatically when the collection of peripheral groups changes, even amongst Kleinian groups, we describe three examples where the group is $F_2$. Each relatively hyperbolic pair $(G, \mathcal{P})$ can be realized as a Kleinian group, where the collection $\mathcal{P}$ is parabolic, but the peripheral subgroups are different, which changes the relatively hyperbolic boundary.  Each has peripheral groups consisting of all the conjugates of some subset of the elements corresponding to the curves $a$, $b$ and $c$ on the one-holed torus:

\begin{figure}[H]
  \centering
  \includegraphics[width=.2\linewidth]{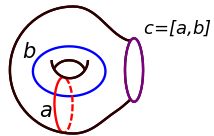}
  \caption{A torus with a boundary component}
  \label{fig:sub1}
\end{figure}

\begin{enumerate}
	\item $G=F_2 = \langle a, b \rangle$, $\P= \varnothing$, $X=$ convex hull of the limit set.

	$\partial(F_2, \varnothing) = \mathcal{C} = $ Cantor set. This is the same as the Gromov boundary, since the set of peripheral subgroups is empty. Here $a, b$ are realized as isometries of $\mathbb{H}^2$ that map the bottom black curve to the top and the left to the right in the central octagon respectively. The group acts geometrically on the convex hull of the limit set (region enclosed by the red axes of conjugates of $\langle[a, b]\rangle$), as shown in Figure \ref{fig:cantor}.

	\begin{figure}[H]
	\centering
	\includegraphics[width=.4\linewidth]{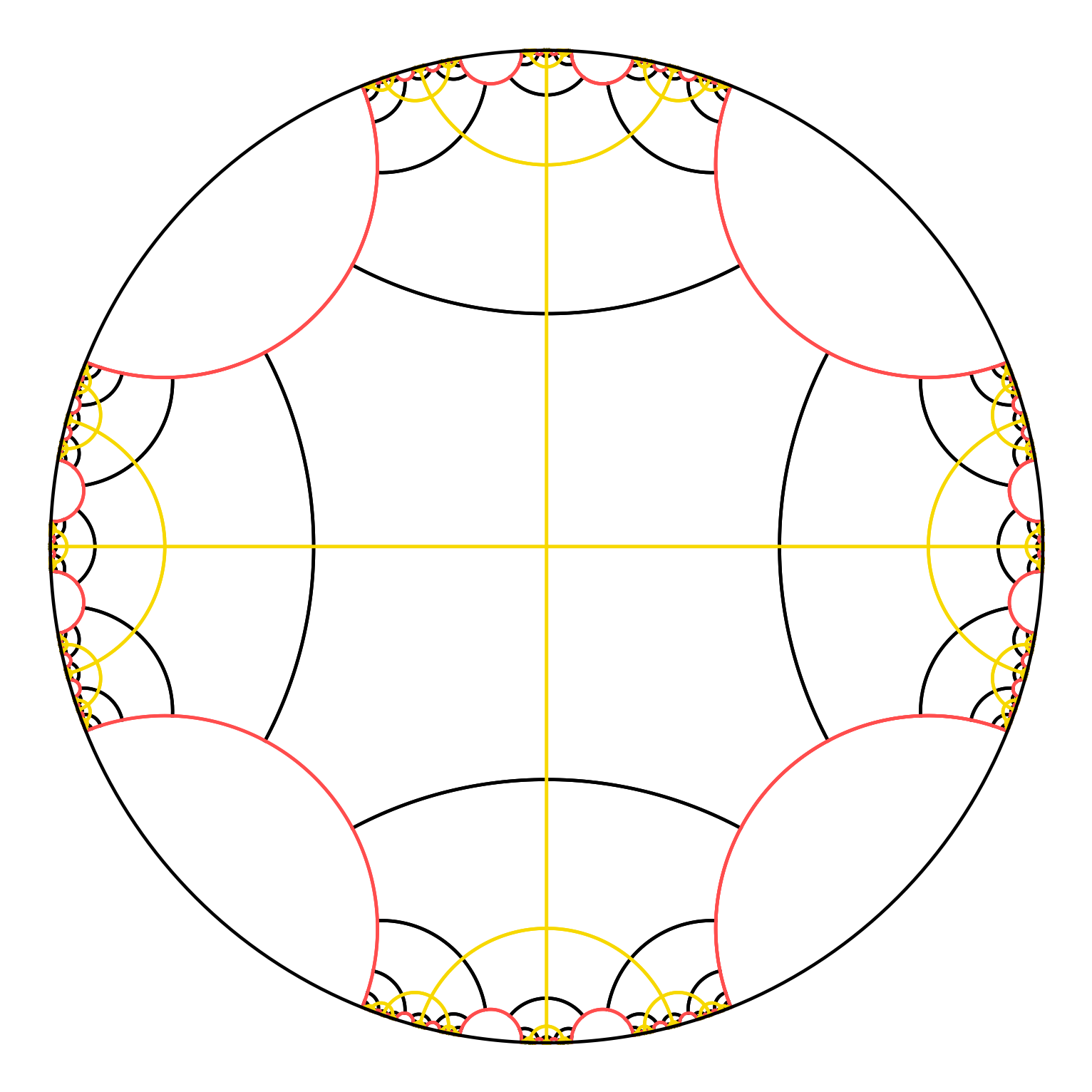}
	\caption{a (relatively) hyperbolic action of $(F_2, \P)$ on $\mathbb{H}^2$ \\ with quotient a torus with a boundary component}
	\label{fig:cantor}
	\end{figure}

	\item (See Figure \ref{fig:ptorusgp}.) $G=F_2=\langle a,b\rangle$, $\P$ is the collection of conjugates of the subgroup $\langle [a,b]\rangle$,  $X=\H^2$.

	$\partial (F_2,\langle [a,b] \rangle) = S^1$.  This can be realized by putting a finite-area hyperbolic structure on the cusped torus.  With such a representation, $F_2$ is a finite co-volume subgroup of $\operatorname{Isom}(\mathbb{H}^2)$ and the limit set and the Bowditch boundary are both $S^1$.

	\begin{figure}[H]
	\centering
	\includegraphics[width=.4\linewidth]{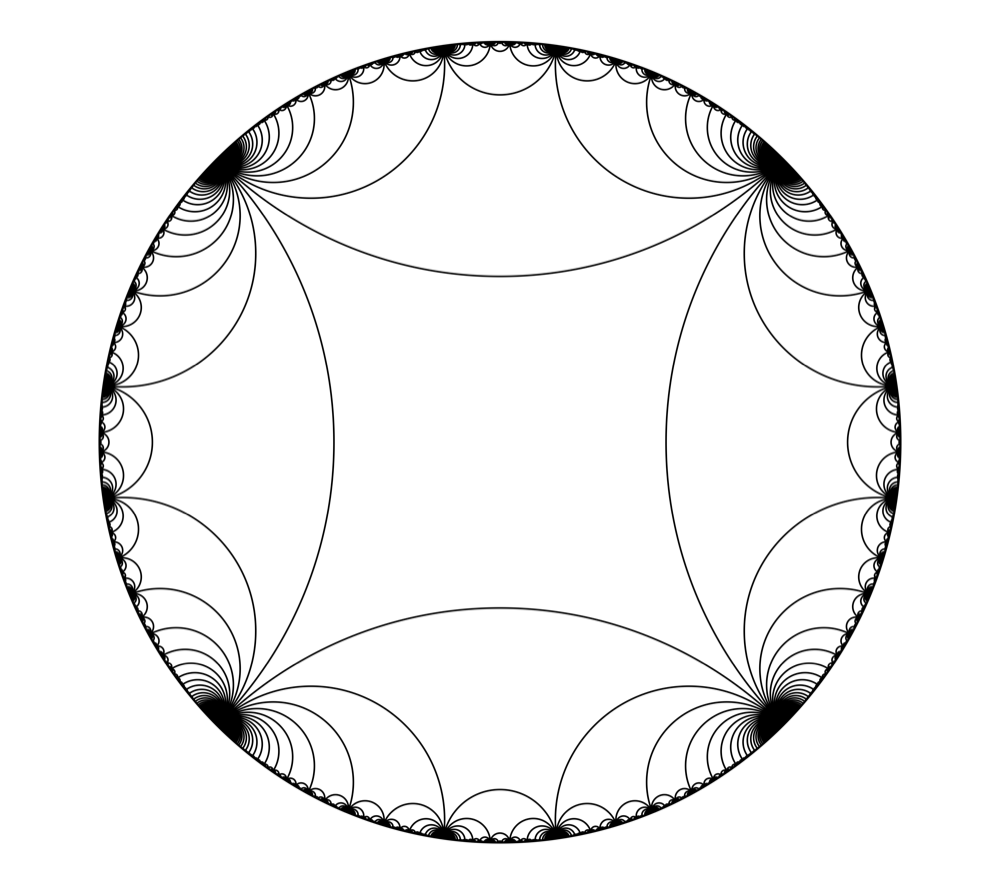}
	\caption{a (relatively) hyperbolic action of $(F_2, \P)$ on $\mathbb{H}^2$ \\ with quotient a cusped torus}
	\label{fig:ptorusgp}
	\end{figure}

	\item (See Figure \ref{fig:gasket}.) $G=F_2$, $\P$ is the collection of all conjugates of the subgroups $ \{ \langle a \rangle, \langle b \rangle, \langle [a,b] \rangle \}$.

	$\partial (G,\P) = $ Apollonian gasket, see \cite{HPWKleinainBoundary}.  Here $X$ is the convex hull of the Apollonian gasket in $\mathbb{H}^3$ and $G$ acts as a geometrically finite Kleinian group on $X$.

			\begin{figure}[h!]
	    \centering
	\includegraphics[width=.8\linewidth]{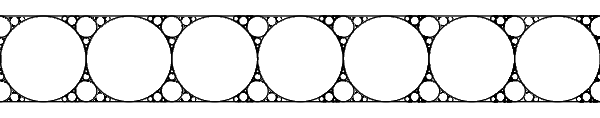}
	    \caption{The Apollonian gasket}
	    \label{fig:gasket}
	\end{figure}

\end{enumerate}
\end{example}

Another source of examples of relatively hyperbolic groups are certain $\cat(0)$ groups.  In particular,  $\cat(0)$ groups with isolated flats admit a relatively hyperbolic group structure.

\begin{defn}
A \textit{flat} is an isometric embedding of $\E^n$ for $n\geq 2$.
\end{defn}

\begin{defn}[Isolated Flats] \label{IsolatedFlats}
        Let $X$ be a $\cat(0)$ space admitting a geometric action by $G$. The space $X$ has \textit{isolated flats} if there is a $G$-invariant collection $\F$ of flats in $X$ such that:

        \begin{enumerate}
            \item There is a constant $D<\infty$ such that each flat in $X$ lies in a $D$-tubular neighborhood of some flat $F\in \F$.
            \item For each positive $\rho<\infty$, there is a constant $\kappa=\kappa(\rho)$ such that for any two distinct flats $F,F'\in \F$, we have

            $$diam(\mathcal{N}_\rho(F) \cap \mathcal{N}_\rho(F') ) < \kappa$$
        \end{enumerate}
\end{defn}
The first condition says that all flats of $X$ are close to a flat of $\F$ and therefore we think of $\F$ as the collection of maximal flats. The second condition says that the flats in $\F$ are far apart, hence isolated.

\begin{thm}[\cite{HK05}]
If $G$ is a $\cat(0)$ group with isolated flats, then $(G,\P)$ is relatively hyperbolic, where $\P$ is the collection of maximal flat stabilizers.

When $G$ is $\cat(0)$ with isolated flats, the visual boundary $\bnd G$ is well-defined.
\end{thm}

A theorem of Tran relates the two boundaries of a $\cat(0)$ group with isolated flats.

\begin{thm}[\cite{Tran13}] \label{thm:TranCAT}
Let $G$ be a $\cat(0)$ group acting geometrically on $X$ with isolated flats and let $\P$ be the collection of maximal flat stabilizers. There is a surjective map

$$\pi \colon \bnd X \to \bnd (G,\P)$$
which is defined by collapsing the boundary of each maximal flat to a single point.
\end{thm}

\begin{example} The fundamental group of a hyperbolic knot complement. Let $G=\pi_1(S^3\setminus \text{figure 8-knot})$. This is not a hyperbolic group since it contains a $\mathbb{Z} \oplus \mathbb{Z}$.  However, $G$ acts geometrically on truncated $\H^3$, since the peripheral subgroups in the Kleinian structure preserve a collection of horoballs. By a result of Ruane \cite{RuaneTruncHyp} the $\cat(0)$ visual boundary $\bnd G$ is $\mathcal{S}$, the Sierpinski carpet. The group $G$ is $\cat(0)$ with isolated flats.  We can apply Tran's theorem above to see that $(G, \mathcal{P})$ where $\mathcal{P}$ is the collection of peripheral $\mathbb{Z} \oplus \mathbb{Z}$ has Bowditch boundary homeomorphic to $S^2$. This follows from the fact that a decomposition which is a null-sequence is an upper semicontinuous decomposition, \cite[page 14]{Daverman} and a theorem of Moore \cite{Mooreusc} that the quotient of an upper semicontinuous decomposition into non-separating continua of $S^2$ is again $S^2$.
\end{example}

Tran's theorem works when $G$ is hyperbolic (see also \cite{Jasonrelhyp} in these notes):

\begin{thm}[\cite{Tran13}] \label{thm:TranHyp}
Let $(G,\P)$ be a relatively hyperbolic group and let $G$ be hyperbolic. Then there is a surjective map

$$\pi:\bnd G \to \partial (G,\P)$$
defined by collapsing all the boundaries of the $P\in\P$.
\end{thm}

\begin{example}
Let $M^3$ be a hyperbolic manifold with totally geodesic boundary. Let $G = \pi_1(M^3)$. Then the Gromov boundary is $\bnd G\cong \mathcal{S}$, the Sierpinski carpet, and $\bnd (G,\P) \cong S^2$.  This can be seen by collapsing the boundaries of the circles removed from the Sierpinski carpet.
\end{example}

The next examples illustrate the use of Tran's theorem to understand the relatively hyperbolic boundaries of relatively hyperbolic group pairs.

\begin{figure}[H]
    \begin{subfigure}{.5\textwidth}
        \centering
	    \includegraphics[scale=.2]{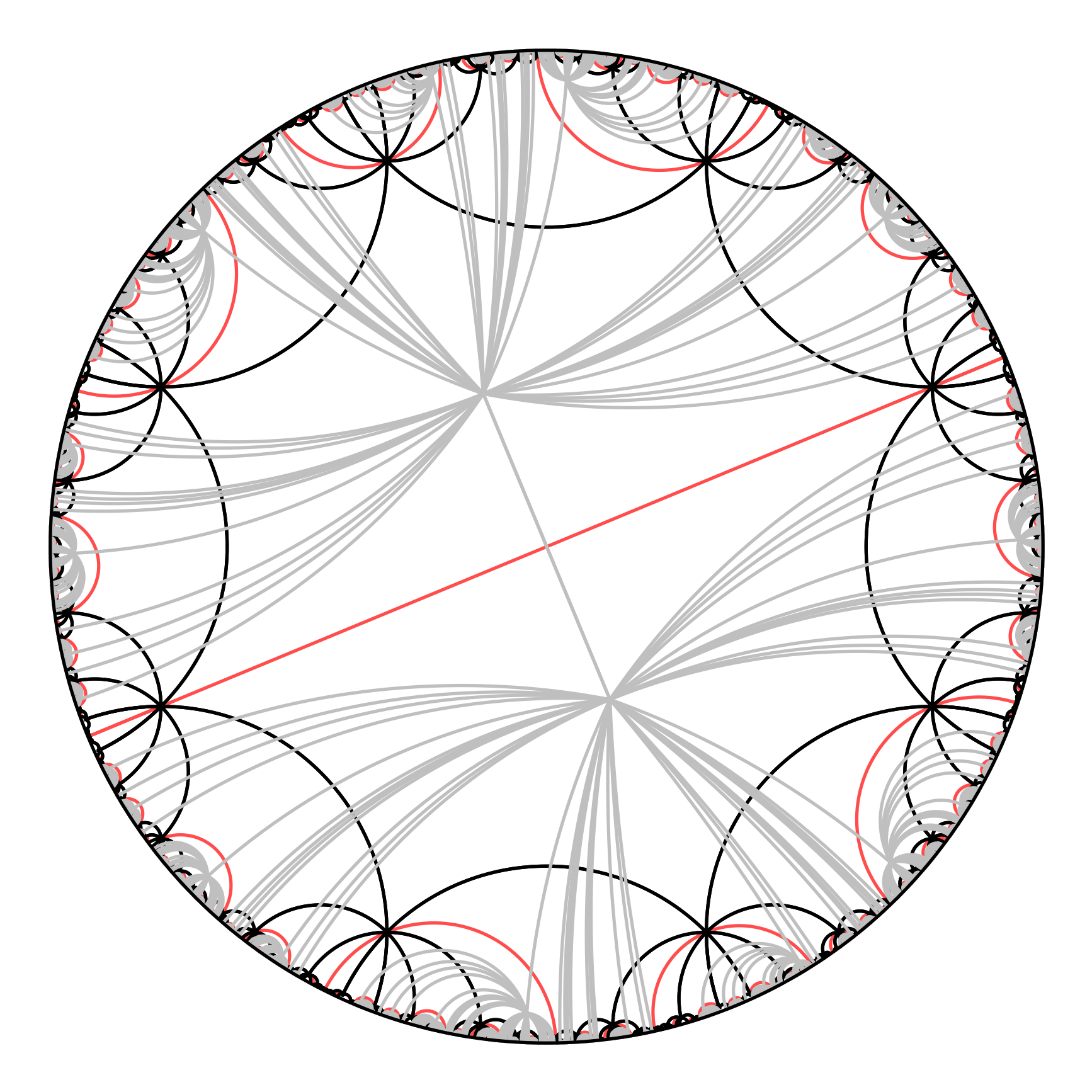}
	    \caption{Bass-Serre tree in $\mathbb{H}^2$}
	    \label{fig:bass-serre}
    \end{subfigure}
    \hspace*{-5pt}
    \begin{subfigure}{.5\textwidth}
	    \centering
	    \includegraphics[width=1.782in]{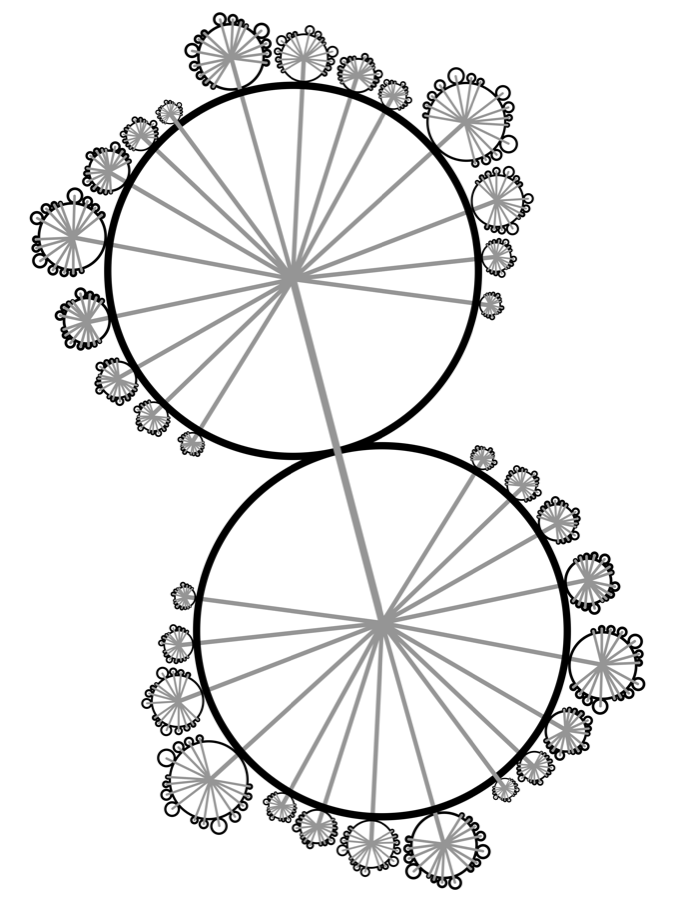}
	    \caption{The Bowditch boundary of $(G, \P)$}
	    \label{fig:tree-of-circles}
    \end{subfigure}
\caption{}

\end{figure}

\begin{example} \label{example:TreeOfCircles} Let $G$ be the fundamental group of a genus 2 surface, realized as a Fuchsian group. We can denote by $c$ the element of $G$ corresponding to a separating curve which bounds a commutator on both sides. Since the subgroup $\langle c \rangle$ is quasi-convex and its conjugates (corresponding to the red separating curves in Figure \ref{fig:bass-serre}) form a malnormal collection, $(G, \mathcal{P})$ is a relatively hyperbolic structure on $G$, where $\mathcal{P}$ consist of $\langle c \rangle$ and its conjugates. With this relatively hyperbolic structure, the Bowditch boundary of $(G,\P)$ is a tree of circles. This boundary $\partial(G, \mathcal{P})$ can be realized by looking at the Bass-Serre tree for the splitting of $G$ over $\langle c\rangle$, where each vertex in the Bass-Serre tree corresponds to a circle in the Bowditch boundary. Two circles meet at a point exactly when there is an edge in the Bass-Serre tree between the vertices. See Figure \ref{fig:tree-of-circles}.
\end{example}

\begin{example}

Here is an example which is not a 3-manifold group.  However, the Bowditch boundary will be planar and we can understand this using Tran's theorem.  Consider three surfaces (with genus at least 1) each with a boundary component. Attach the three boundary curves to the curves of the $T^2$ pictured in Figure \ref{fig:torus 3 }. Let $G$ be the fundamental group of this 2-complex and $\P$ the collection of abelian subgroups of rank 2. The resulting Bowditch boundary is a tree of circles and is planar. Work of Hruska-Walsh \cite{HW1} shows that the $\cat(0)$ boundary contains a $K_{3,3}$ graph. Such a graph is an obstruction to the group acting properly discontinuously on a contractible 3-manifold by work of  Bestvina-Kapovich-Kleiner \cite{BestvinaKapovichKleiner02}.

Using Tran's Theorem \ref{thm:TranCAT}, $\bnd (G,\P)$ is obtained by collapsing the circles in the boundary coming from the $\mathbb{Z} \oplus \mathbb{Z}$ subgroups. This boundary is planar, but has cut points.  See \cite{HW1} for details.

\begin{figure}[H]
    \centering
    \includegraphics[scale=.3]{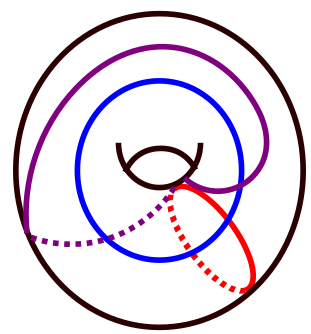}
    \caption{3 curves on $T^2$}
    \label{fig:torus 3 }
\end{figure}
\end{example}

\section{More definitions of relatively hyperbolic groups}\label{day2}

In this section we introduce several more definitions of a relatively hyperbolic group pair. By work of Dahmani, Hruska, and Groves--Manning, all these are equivalent (and are equivalent to our first definition in Section \ref{day1}) \cite{DahmaniRelHypEquiv, HruskaRelHypEquiv, GMCusp}. This gives us multiple ways of identifying and studying relatively hyperbolic groups. Furthermore, two of the definitions are constructive in the sense that algebraic information about the group can be used to build an appropriate hyperbolic space.

As noted before, all groups we consider are finitely generated throughout this section. Many of these definitions can be adapted to non-finitely generated groups, but we do not do so here. See \cite{HruskaRelHypEquiv} for some of these definitions. We recall the definition of a Cayley graph.

\begin{defn}
        For a group $G$ with a generating set $\A$, the \textit{Cayley Graph}, denoted $C(G,\A)$ is a graph with
        \begin{enumerate}
            \item a vertex for each $g\in G$ and
            \item an edge labelled by $a\in \A$ joining the vertices $g$ and $ga$.
        \end{enumerate}
\end{defn}

\begin{figure}[H]
    \centering
\begin{tikzpicture}
\draw (0,0) --(1,1);
\node [left] at (.6,.6) {$a$};
\node [left] at (0,0) {$g$};
\node [right] at (1,1) {$ga$};
\draw[fill] (1,1) circle [radius=0.025];
\draw[fill] (0,0) circle [radius=0.025];
\draw (4,0) -- (5,1);
\node [left] at (4,0) {$bg$};
\node [right] at (5,1) {$bga$};
\draw[fill] (5,1) circle [radius=0.025];
\draw[fill] (4,0) circle [radius=0.025];
\node [left] at (4.6,.6) {$a$};
\draw [|->] (2,.5) -- (3,.5);
\node [above] at (2.5,.5) {$b$};

\end{tikzpicture}
    \caption{The action of $G$ on its Cayley graph}
    \label{fig:cayey-action}
\end{figure}
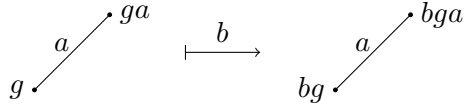

The group $G$ acts on its Cayley graph on the left as seen in Figure \ref{fig:cayey-action}

\vskip .2in

For hyperbolic groups, the Cayley graph, endowed with a metric where each edge has length 1, captures the hyperbolic geometry. For a relatively hyperbolic group, however, the Cayley graph does a poor job of capturing desired geometric properties of the group. Farb introduced the notion of the coned-off Cayley graph as a way of capturing these properties using a graph akin to the Cayley graph.

\begin{defn}[Coned-off Cayley Graph, \cite{FarbCoset}]
        Let $(G,\P)$ be a relatively hyperbolic pair and $\A$ a finite generating set for $G$. The \textit{coned-off Cayley graph}, denoted $\mathcal{C}(G,\P,\A)$ is the Cayley graph $C(G,\A)$ with some additions. Let $\mathbb{P}$ be a set of representatives for the conjugacy classes of $\P$. For each coset $gP$, where $P\in \mathbb{P}$, add a vertex $v_{gP}$. Then for each $h\in gP$, add an edge of length $\frac{1}{2}$ from $h$ to $v_{gP}$. The resulting graph is $\mathcal{C}(G,\P,\A)$. If $\gamma$ is a path in $C(G,\A)$, then let $\hat \gamma$ be the path in $\mathcal{C}(G,\P,\A)$ where we replace each maximal subpath in a coset of a peripheral subgroup with two edges of length 1/2, meeting $v_{gP}$.
\end{defn}

In an ideal world, we would say when this resulting graph is hyperbolic, then the group is relatively hyperbolic. Unfortunately, this would be too broad of a definition, as it would allow $\Z^2$ to be relatively hyperbolic, as shown in the next example.

\begin{example}\label{example:Z rel hyp}
Consider $G=\langle a,b| [a,b]\rangle = \Z\oplus \Z$ and let $\P$ consist of $\langle a \rangle$ and its conjugates (which is just $\langle a \rangle$). Then the coned off Cayley graph $\mathcal{C}(G,\P,\A)$ is hyperbolic (as it is quasi-isometric to a line) where $\A = \{ a, b\}$.
\end{example}

In light of this example, the following definitions are required to achieve the desired definition.

\begin{defn}[Without backtracking]
        A path in $\mathcal{C}(G,\P,\A)$ is \textit{without backtracking} if once the path hits $v_{gP}$, it never returns. If $\gamma$ is a path in $C(G, \mathcal{A})$, we say it is without backtracking if $\hat\gamma$ is without backtracking. A path in $\gamma$ in $C(G, \mathcal{A})$ \textit{penetrates} the coset $gH$ if $\hat\gamma$ passes through $v_{gH}$.
\end{defn}

\begin{defn}[Bounded coset penetration]
        Let $C(G,\P,\A)$ be the coned-off Cayley graph for $(G,\A)$. This has \textit{bounded coset penetration} if for each $\lambda\geq 1$, there is a constant $a(\lambda)>0$ such that if $\gamma,\gamma'$ are two $(\lambda,0)$-quasigeodesics in $C(G,\A)$ without backtracking, with the same initial vertex, and with endpoints that are no more than $1$ apart, then the following two conditions hold:

        \begin{enumerate}
            \item if $\gamma$ penetrates $gP$ and $\gamma'$ does not, then the entering and exiting vertices (endpoints of the subpath in $gP$) of $\gamma$ are at most $a(\lambda)$ from each other in $C(G,\A)$.
            \item if $\gamma$ and $\gamma'$ both penetrate $gP$, then, in $C(G,\A)$, the two entering vertices of each are $a(\lambda)$-close and the two exiting vertices are $a(\lambda)$-close.
        \end{enumerate}
\end{defn}
Example \ref{example:Z rel hyp} does not satisfy the bounded coset penetration, which is exactly what we wanted. Therefore, the following is one of our definitions of relatively hyperbolic:

\begin{defn}[Relatively hyperbolic]\label{def: coned-off}
        A group pair $(G,\P)$ is relatively hyperbolic if the coned-off Cayley graph $\mathcal{C}(G,\P,\A)$ is $\delta$-hyperbolic and satisfies the bounded coset penetration property \cite{FarbCoset}.
\end{defn}

The equivalence of Definition \ref{def: coned-off} with Definition \ref{def:original rel hyp} was proved by Hruska \cite{HruskaRelHypEquiv}.

Something to note: the hyperbolic space $\mathcal{C}(G,\P,\A)$ is \textit{not} proper because the vertices $v_{gP}$ have infinite valence. We can still define the boundary of $\mathcal{C}$ as before, but it is not compact because $\mathcal{C}$ is not proper. This boundary is missing the parabolic points.  There is, however, a relationship between the Bowditch boundary and the boundary of the coned-off Cayley graph, see \cite[Theorem 9.1]{Bowrelhyp} for the correct topology:

$$\partial (G,\P) = \partial (\mathcal{C}(G,\P,\A)) \cup \bigcup_{P \in \mathbb{P}, g\in G} \{v_{gP}\}$$

\begin{example}
Let $G=\pi_1(S^2\setminus K)$ be a hyperbolic knot complement. Then $\partial (G, \mathcal{P}) = S^2$ because $G$ is a geometrically finite Kleinian group with finite co-volume. The set of parabolic fixed points is dense in $S^2$.  So we can understand the boundary of the coned-off Cayley graph $  \partial (\mathcal{C}(G,\P,\A))$ as $S^2$ with a countable dense collection of points removed. These are the parabolic fixed points.
\end{example}

The next definition of a relatively hyperbolic pair also uses the Cayley graph to construct an appropriate hyperbolic space. But first we introduce combinatorial horoballs:

\begin{defn}[Combinatorial Horoballs, \cite{GMCusp}]
Let $\Gamma$ be a graph with all edges length 1. Construct a new graph $\mathcal{H}(\Gamma)$ with vertex set
$$V(\mathcal{H}(\Gamma)) = V(\Gamma) \times \Z_{\geq 0}.$$
There are two types of edges in $\mathcal{H}(\Gamma)$: For all $v\in V(\Gamma)$, there is an edge between $(v,k)$ and $(v,k+1)$. For each $k$, there is an edge between $(v,k)$ and $(w,k)$ if, in $\Gamma$, $d(v,w)\leq 2^k$. The first type of edges we call vertical and the second type are horizontal.

\end{defn}
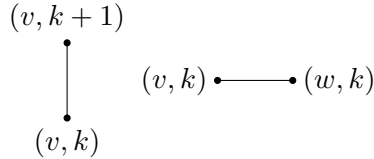
\begin{figure}[H]
    \centering
\begin{tikzpicture}

    \node [below] at (0,0) {$(v,k)$};
    \node [above] at (0,1) {$(v,k+1)$};
    \draw (0,0) -- (0,1);
    \draw[fill] (0,0) circle [radius=0.04];
    \draw[fill] (0,1) circle [radius=0.04];

    \node[left] at (2,.5) {$(v,k)$};
    \node[right] at (3,.5) {$(w,k)$};
    \draw (2,.5) -- (3,.5);
    \draw[fill] (2,.5) circle [radius=0.04];
    \draw[fill] (3,.5) circle [radius=0.04];

    \end{tikzpicture}
    \caption{Vertical and horizontal edges in the graph $\mathcal{H}(\Gamma)$}
    \label{fig:horoballs}
\end{figure}

The motivation for this definition comes from the example of $F_2$ acting on $\H^2$ with quotient a cusped torus.  In this action, there is a collection of invariant horoballs centered at the parabolic points on the boundary. The combinatorial horoball definition and the original Gromov definition, see \cite{Gromov87, Szcz} model this behavior.

Given a group pair $(G, \P)$, we can construct a graph using the Cayley graph and combinatorial horoballs. If the resulting space is hyperbolic, then the group pair $(G,\P)$ is relatively hyperbolic, and the combinatorial horoballs mimic the behavior seen above in $F_2$ acting on $\H^2$.

Note, in the definition below, $\mathbb{P}$ is a \textbf{finite} collection of parabolic subgroups, which differs from the definition of $(G,\P)$ from Section \ref{day1}. To go from $\mathbb{P}$ to $\P$, take the union of all conjugacy classes of $P\in \mathbb{P}$. To go the other way, pick a representative of each conjugacy class in $\P$.

\begin{defn}[Cusped Cayley graph, \cite{GMCusp}]\label{def:cusped-off}
Let $G$ be a group and $\mathbb{P}$ a finite collection of subgroups. Let $\A$ be a generating set for $G$ which contains a generating set for each $P\in \mathbb{P}$. Construct the Cayley Graph $C(G,\A)$ and, for each coset $gP$ of each $P\in \mathbb{P}$, attach a copy of $\mathcal{H}(gP)$ to $C(G,\A)$. Here the $0$-level of $\mathcal{H}(gP)$ is identified with $C(gP)$. We denote this space $X(G,\mathbb{P})$ and call it the \textit{cusped Cayley graph}.
\end{defn}

Note that this construction requires a generating set for \textbf{both} $G$ and each $P\in \mathbb{P}$, so each parabolic subgroup must be finitely generated. However, if $(G, \P)$ is relatively hyperbolic and $G$ is finitely presented, then each $P\in \mathbb{P}$ is finitely presented as well \cite{DahmaniPresenting} .

\begin{thm}[Groves-Manning, \cite{GMCusp}]
The pair $(G,\P)$ is relatively hyperbolic (in the sense of Definition \ref{def:original rel hyp}) when the cusped Cayley graph $X(G,\mathbb{P})$ is hyperbolic. Furthermore, $\partial (X(G,\mathbb{P}))$ is the Bowditch boundary.
\end{thm}

Both Definitions \ref{def: coned-off} and \ref{def:cusped-off} give constructions for creating a hyperbolic space admitting an action by the relatively hyperbolic group pair $(G,\P)$. By taking algebraic information (a group, a collection of subgroups, and a generating set), we can build a geometric model for $(G,\P)$, which is not the case for Definition \ref{def:original rel hyp}. Furthermore, the boundaries of the resulting spaces are either very close to the Bowditch boundary (in the coned-off Cayley graph) or exactly the Bowditch boundary (in the cusped Cayley graph). And lastly, unlike the spaces satisfying Definition \ref{def:original rel hyp}, any two cusped Cayley graphs for the same relatively hyperbolic pair are quasi-isometric \cite{HruskaHealyCusped}.

Our final definition of relatively hyperbolic is also due to Bowditch (as was the original). As in Definitions \ref{def: coned-off} and \ref{def:cusped-off}, the hyperbolic space of interest will be a graph. But unlike those definitions, the graph does not come from the Cayley graph, nor is it constructive.

\begin{defn}
        A graph $K$ is \textit{fine} if each edge of $K$ is contained in only finitely many circuits of length $n$ for each $n$.
\end{defn}

\begin{defn}[Relatively hyperbolic, \cite{Bowrelhyp}]\label{def:fine rel hyp}
Let $G$ act on a $\delta$-hyperbolic graph $K$ with finite edge stabilizers and finitely many orbits of edges. If $K$ is fine, then $(G,\P)$ is \textit{relatively hyperbolic}, where $\P$ consists of stabilizers of infinite valance vertices.
\end{defn}

The equivalence of this definition with Definition \ref{def:original rel hyp} is due (independently) to Bowditch \cite[Theorem 7.10]{Bowrelhyp}, Dahmani \cite{DahmaniRelHypEquiv}and Hruska \cite{HruskaRelHypEquiv}.

We have 4 equivalent, yet separate, definitions of relatively hyperbolic, which we will summarize here.  Note there are more (equivalent) definitions, for example \cite{gerasimov, Yaman}.

\begin{enumerate}
    \item Definition \ref{def:original rel hyp}: When $G$ acts properly discontinuously and by isometries on a hyperbolic metric space $X$ with each $x\in \bnd X$ either conical limit point or bounded parabolic and $\P$ is the collection of maximal parabolic subgroups, then $(G,\P)$ is relatively hyperbolic.
    \item Definition \ref{def:fine rel hyp}: If $G$ acts on a $\delta$-hyperbolic graph $K$ with finite edge stabilizers and finitely many orbits of edges, and $K$ is fine, then $(G,\P)$ is relatively hyperbolic, where $\P$ is the collection of stabilizers of infinite valance vertices.

    \item Definition \ref{def: coned-off}: When the coned-off Cayley graph $\mathcal{C}(G,\P,\A)$ is hyperbolic and has bounded coset penetration, $(G,\P)$ is relatively hyperbolic.
    \item Definition \ref{def:cusped-off}: When the cusped Cayley graph $X(G,\mathbb{P})$ is hyperbolic, $(G,\P)$ is relatively hyperbolic.

\end{enumerate}
\section{What the boundary tells us and the relation to Kleinian groups}\label{day3}

What can the boundary of a relatively hyperbolic group tell you about the group?

We'll begin by examining the case of Kleinian groups.  These are key examples of hyperbolic and relatively hyperbolic groups.
Many of the results about hyperbolic and relatively hyperbolic groups in this section were inspired by known results regarding Kleinian groups and their associated manifolds.  We will be discussing some manifold theory without always giving detailed definitions. The first chapter of \cite{Kapbook} has comprehensive definitions.
The key idea is that a hyperbolic three-manifold has a ``characteristic submanifold" containing the essential annuli in the manifold.  This can be seen from the limit set of the associated Kleinian group.  A similar and important phenomenon happens with hyperbolic and relatively hyperbolic groups.  This theory was begun by Bowditch (echoing the Jaco-Shalen and Johannsen characteristic submanifold theory) and continued by many people.

\begin{defn}[Kleinian group]
A group $\Gamma$ is \textit{Kleinian} if it is a discrete subgroup of $\PSL(2,\mathbb{C})$. Note that $\PSL(2,\mathbb{C}) = \operatorname{Isom}^+(\mathbb{H}^3)$, the orientation preserving isometries of $\mathbb{H}^3$.
\end{defn}

\begin{defn}[Limit set]
Let $\Gamma$ be a Kleinian group and let the boundary of $\H^3$ be $\hat{\C}$. Fix $x\in \H^3$, then the \textit{limit set} of $\Gamma$, denoted $\Lambda_\Gamma$, is

$$\Lambda_\Gamma:=\overline{\Gamma \cdot x} \cap \hat{\C}$$

\end{defn}

\begin{remark}
    The choice of $x\in \H^3$ does not change $\Lambda_\Gamma$.
\end{remark}

\begin{defn}[Geometrically finite]
Let $\Gamma$ be a Kleinian group and let $C(\Gamma)$ be the convex hull of $\Lambda_\Gamma$. Then $\Gamma$ is \textit{geometrically finite} if $C(\Gamma)/\Gamma$ has finite volume.

\end{defn}

A natural connection between Kleinian groups and relatively hyperbolic groups comes from the following fact: if $\Gamma<\PSL(2,\mathbb{C})$ is geometrically finite, then $(\Gamma, \P)$ is relatively hyperbolic and $\Lambda_\Gamma = \bnd (G,\P)$, where $\P$ is the collection of parabolic subgroups of $\Gamma$, see \cite{BowGF}.

An example and a non-example of geometrically finite Kleinian groups:

\begin{example} If $\Gamma < \PSL(2, \mathbb{C})$ is $\pi_1(M^3)$ where $M^3$ is a closed hyperbolic manifold, then $\Gamma$ is geometrically finite. Also, when $\H^3/\Gamma$ is a finite-volume hyperbolic cusped 3-manifold, $\Gamma$ is geometrically finite.  More generally, if $H< \Gamma$ is a finite index subgroup, where $\Gamma$ is a geometrically finite Kleinian group, then $H$ is also geometrically finite. Note that when $H$ is finite index  then  $\Lambda_H=\Lambda_\Gamma$, and $H$ acts geometrically finitely on the convex hull of its limit set.
\end{example}

\begin{example}
Non-Example: Let $\psi$ be a pseudo-Anosov homeomorphism of a hyperbolic surface $S_g$.  The mapping torus $M_\psi^3$ is a hyperbolic 3-manifold and the limit set of its fundamental group, $\Gamma$, is $\hat{\bbC}$. This manifold fibers over the circle and $\pi_1(S_g)$ is normal in $\Gamma$, so $\Lambda_{\pi_1(S_g)}=\Lambda_\Gamma=\hat{\bbC}$. Since $\Lambda_\Gamma$ is the entire boundary of $\H^3$, the convex hull of the limit set is all of $\H^3$.  $\H^3/\pi_1(S_g)$ is infinite volume, hence this action is not geometrically finite.
\end{example}

For some time, lots of manifolds have been known to admit geometrically finite hyperbolic structures.  By work of Thurston, Haken manifolds whose fundamental groups do not contain free abelian groups of rank 2 can be realized as hyperbolic manifolds. The proof is quite involved, see \cite{Kapbook} and \cite{MorganOnUniformization}, and includes the case when the manifold has boundary. A 3-manifold $M$  is \emph{irreducible} if every 2-sphere in $M$  bounds a 3-ball in $M$ and \emph{atoroidal} if it is irreducible and $\pi_1(M)$ contains no $\mathbb{Z} \oplus \mathbb{Z}$.


\begin{thm}[Theorem A, page 70 \cite{MorganOnUniformization}]
Let $M$ be a compact, atoroidal, Haken 3-manifold. Then there is a geometrically finite, complete hyperbolic manifold $N$ such that $C(\pi_1(N))/\pi_1(N)$ is homeomorphic to $M$.
\end{thm}

There is a corresponding theorem for ``pared manifolds" \cite[pg70 Theorem B']{MorganOnUniformization}].  This theorem provides a large family of relatively hyperbolic pairs, where the peripheral subgroups are exactly the parabolic subgroups of the corresponding Kleinian group.

\begin{remark} \label{remark genus 2}
    It is possible for the manifold to have annuli. For example, let $S_{2,1}$ be a  genus two surface with one boundary component.  Then the three manifold obtained by gluing three copies of $S_{2,1} \times I$ along $\partial(S_{2,1}) \times I$ to a solid torus along three parallel annuli on the boundary of the solid torus satisfies the hypotheses above. Thus this admits a geometrically finite hyperbolic structure.  The quotient of $H^3$ by this Kleinian group is infinite volume, however the quotient of the convex hull of the limit set has finite volume.  Thus it acts geometrically finitely on the convex hull of its limit set, which is hyperbolic.
\end{remark}
The following definition allows us to understand the essential annuli in a 3-manifold.  For geometrically finite hyperbolic 3-manifolds, the characteristic submanifold can be described from the limit set, see \cite{walshbump}.

\begin{defn}[Characteristic Submanifold]
Let $(M,\partial M)$ be a 3-manifold with incompressible boundary (e.g. $\pi_1(M)$ is a geometrically finite Kleinian group). Then the \textit{characteristic submanifold}, $(X,S)$, is a submanifold with:

\begin{enumerate}
    \item Each component is an $I$-bundle over a surface or a solid torus with a Seifert-fibered structure.
    \item $\partial X \cap \partial M = S$
    \item The components of $\partial X \setminus S$ are essential annuli.
    \item Any essential annulus or M\"obius band is properly homotopic into $(X,S)$.
    \item $(X,S)$ is unique up to isotopy.

\end{enumerate}

\end{defn}

The characteristic submanifold allows us to detect a splitting of the fundamental group of the 3-manifold (with incompressible boundary) along infinite cyclic subgroups coming from the essential annuli. A result of Bowditch tells us that we can detect such a splitting in \textbf{any} hyperbolic group, and this splitting can be detected through the topology of the boundary \cite{BowHypeSplitting}.

\begin{defn}[Splitting]
A {\it splitting} of a group $G$ over a class of subgroups is a non-trivial finite graph of groups representation of $G$, where each edge group belongs to the class.
\end{defn}

\begin{thm}[\cite{BowHypeSplitting}]\label{thm:BowHyp}
Let $G$ be a hyperbolic group. If $\partial G$ has a local cut point, then $G$ splits over a virtually cyclic (i.e. 2-ended) subgroup. Furthermore, $G$ splits as a bipartite graph of groups with three types of vertices:

\begin{enumerate}
    \item virtually cyclic
    \item virtually Fuchsian
    \item rigid--these contain no further splittings.
\end{enumerate}
\end{thm}

To see this in the boundary $\partial G$ of $G$, there is a cut pair in $\partial G$ which is the limit set of the subgroup that $G$ splits over.  In fact, all the conjugates of this subgroup will have limit set consisting of a cut pair.

\begin{example}
Let $A, B, C$ be three copies of $S \times I$, where $S$ is a torus with one boundary component.  Let $T$ be a solid torus $T^2 \times I$. Glue the $\partial(S) \times I$ of $A, B, C$ to parallel longitudinal annuli on $T \times \{0\}$ by degree 1 maps.  We obtain a 3-manifold $M$ with boundary as shown in Figure \ref{fig:surfacethickened}, where $T$ is tricolor.

The fundamental group $G$ of $M$ has a graph of groups decomposition with three vertex groups $F_2 = \pi_1(A) = \pi_1(B) = \pi_1(C)$, one vertex group $\mathbb{Z} = \pi_1(T)$ and three edge groups $\mathbb{Z}$ (Figure \ref{fig:surfacegraphofgroups}).

\begin{figure}[H]
    \begin{subfigure}{.5\textwidth}
    \hspace*{1cm}
	    \includegraphics[scale=.075]{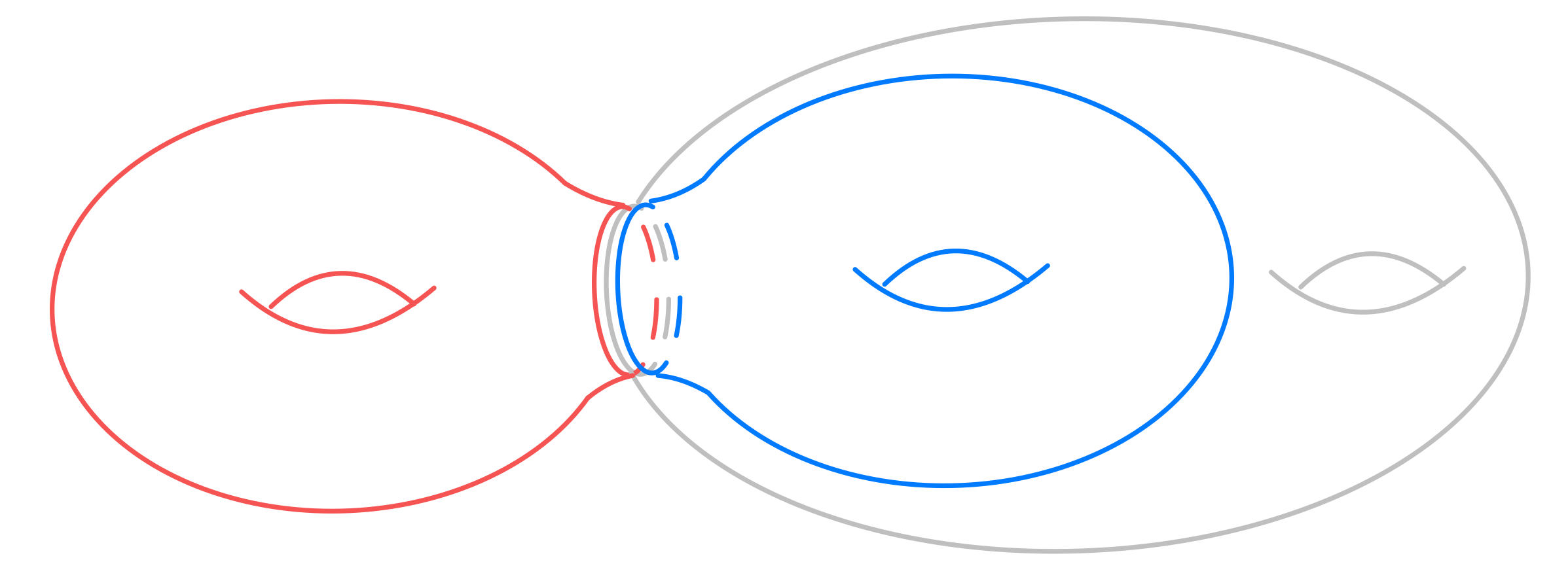}
	    \caption{$M$}
	    \label{fig:surfacethickened}
    \end{subfigure}
    \hspace*{-5pt}
    \begin{subfigure}{.5\textwidth}
	    \centering
	    \includegraphics[scale=.12]{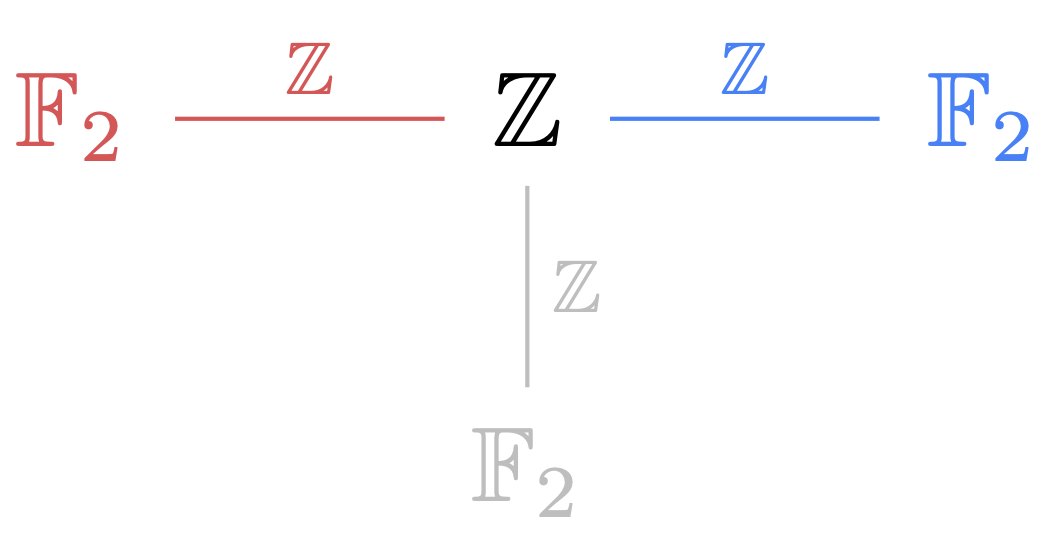}
	    \caption{graph of group decomposition }
	    \label{fig:surfacegraphofgroups}
    \end{subfigure}
\caption{}
\label{fig:surface}
\end{figure}

The universal cover $\tilde M$ of $M$ is a tree of spaces, where the tree is bipartite with two types of vertices. Vertices of type I are universal covers of $T$, and have valence 3. Vertices of type II are universal covers of $A, B, C$, as in Figure \ref{fig:cantor}, and have valence $\infty$. Figure \ref{fig:surfaceuniversalcover} shows three sheets of universal covers of $A, B, C$ attached along a universal cover of $T$, the vertical thickened line. Each sheet contains infinitely many universal covers of $T$, whose fundamental groups are conjugates of $\pi_1(T)$. Attached to each universal cover of $T$ are three sheets of universal covers of $A, B, C$.

For the relatively hyperbolic boundary of $G$, note that $G$ acts geometrically on two proper hyperbolic metric spaces: the convex hull of $\Lambda(G)$ and $\tilde M$, so $\partial G \cong \Lambda(G) \cong \partial (\tilde M)$ (Figure \ref{fig:surfacelimitset}).

\begin{figure}[H]
    \begin{subfigure}{.5\textwidth}
	    \centering
	    \includegraphics[scale=.15]{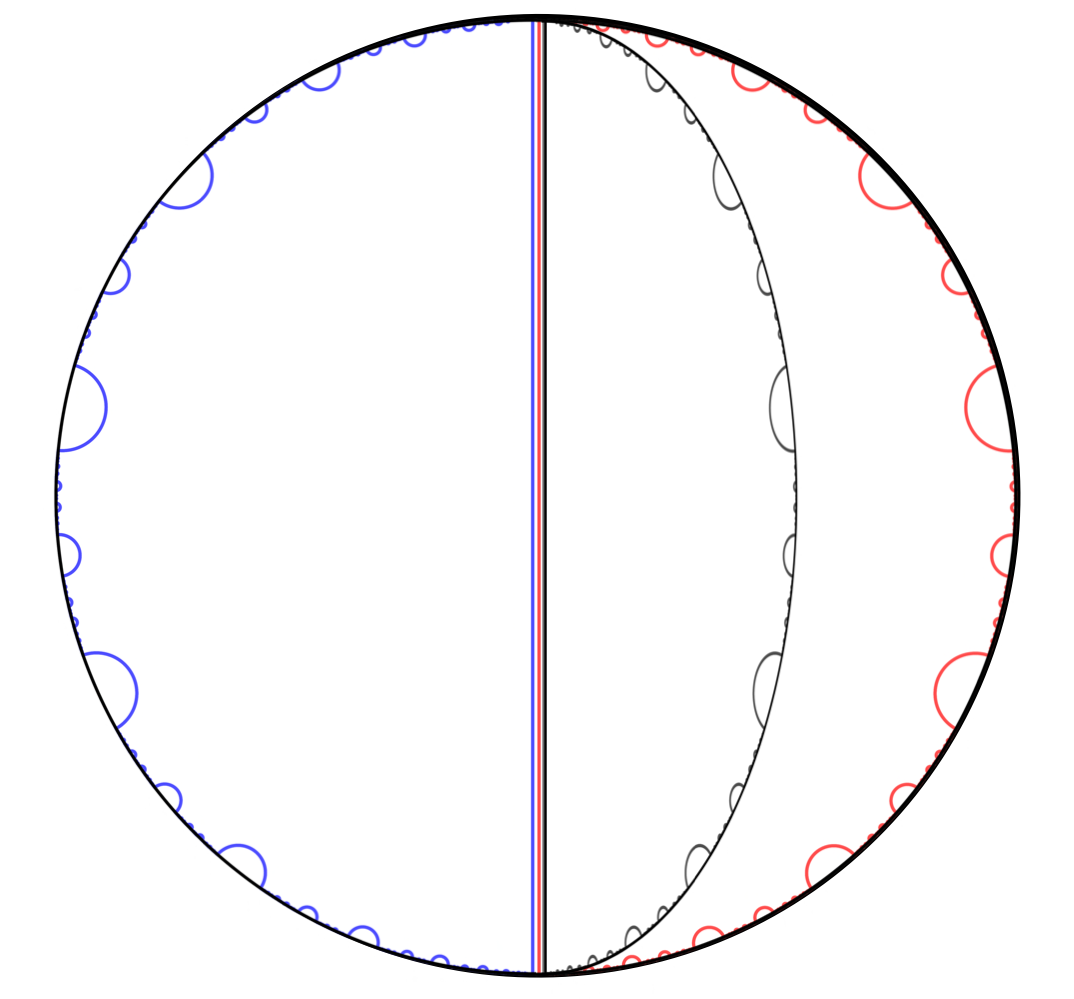}
	    \caption{the universal cover $\tilde M$ of $M$}
	    \label{fig:surfaceuniversalcover}
    \end{subfigure}
    \hspace*{-5pt}
    \begin{subfigure}{.5\textwidth}
	    \centering
	    \includegraphics[width=2.038in]{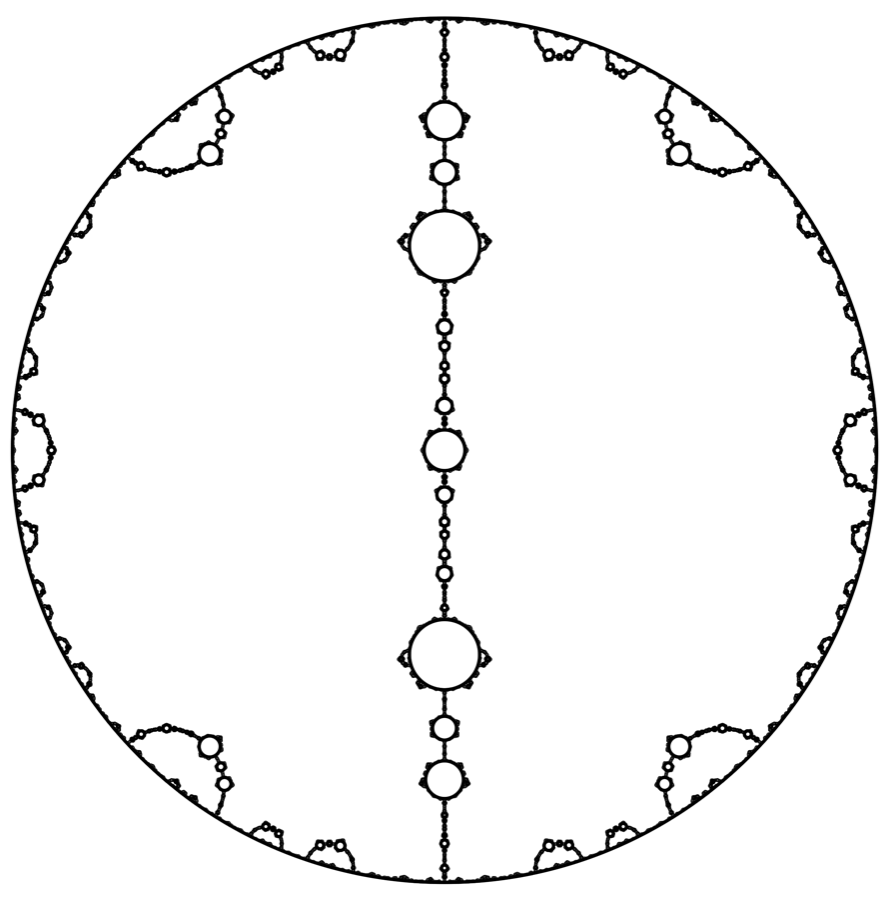}
	    \caption{the limit set $\Lambda(G) \cong \partial(\tilde M)$}
	    \label{fig:surfacelimitset}

    \end{subfigure}
    \caption{}
\end{figure}

In Bowditch's language (Theorem \ref{thm:BowHyp}), the fundamental groups of the 3-valent vertices are of type 1 (virtually cyclic), and those of the $\infty$-valent vertices are of type 2 (virtually Fuchsian). All are quasiconvex subgroups of $G$, with their boundaries embedded in $\partial G$. Here, $\partial G$ has a tree-like structure, where each vertex is a pair of points (for type 1 vertices, of valence 3) or a Cantor set (for type 2 vertices, of valence $\infty$), and where each edge is a pair of points (for the edge groups $\mathbb{Z}$) that coincides with vertices of type 1. In $\Lambda (G) = \partial G$, we see a cut pair of valence 3 at the north and south pole, corresponding to the boundary of the 2-ended subgroup $\pi_1(T) = \mathbb{Z}$ of type 1. All other cut pairs of valence 3 correspond to conjugates of the 2-ended subgroup. We can see $\partial G$ as three Cantor sets glued together along every 3-valent cut pair, together with the boundary of the bipartite tree, which corresponds to rays that keep switching sheets and thus does not belong to the boundary of any vertex or edge.
\end{example}

Kapovich and Kleiner use Bowditch's result to classify the types of 1-dimensional boundaries possible for 1-ended hyperbolic groups.

\begin{thm}[\cite{KapKleinerOneDimBoundary}]\label{thm:KapKleiner}
Let $G$ be a 1-ended hyperbolic group with $\partial G$ 1-dimensional and $G$ does not split over a 2-ended subgroup. Then $\partial G$ is homeomorphic to one of the following:

\begin{enumerate}
    \item $S^1$
    \item $\mathcal{S}$, a Sierpinski carpet (planar)
    \item a Menger curve (non-planar)
\end{enumerate}
\end{thm}
The reason that we have the two-ended hypothesis is that there are Kleinian groups with boundaries different from above (as in Figure \ref{fig:surfacelimitset}).  A hyperbolic manifold group that splits over an infinite cyclic group is not rigid in the sense that it admits many different hyperbolic structures.
Furthermore, work of Canary and McCullough shows that there are hyperbolic 3-manifolds $M_1$ and $M_2$ with $\pi_1(M_1)\cong \pi_1(M_2)$ but $M_1$ and $M_2$ are not homeomorphic.  See \cite{CanaryMcMullough}.  For example, if we alter the example of Remark \ref{remark genus 2} so that the surfaces have different genera, then changing the cyclic order around the solid torus will not change the fundamental group.

\begin{defn}[Peripheral splitting]
Let $(G, \P)$ be a relatively hyperbolic pair. A splitting is {\it relative} to $\P$ if each subgroup in the collection $\P$ is conjugate into one of the vertex groups.

A {\it peripheral splitting} of $(G, \P)$ is a bipartite splitting of $G$ relative to $\P$, where each $P \in \P$ is conjugate into vertices of one color.
\end{defn}

\begin{defn}[Tame]
Let $(G,\P)$ be a relatively hyperbolic pair. A subgroup $P\in \P$ is \textit{tame} if $P$ is finitely generated, 1- or 2-ended, and does not contain an infinite torsion subgroup.
\end{defn}

\begin{defn}[Rigid, \cite{DahmaniRecognizing}]
A relatively hyperbolic pair $(G, \P)$ is {\it rigid} if $G$ has no splitting relative to $\P$ over virtually cyclic groups or over parabolic subgroups.
\end{defn}

Matt Haulmark proves a theorem similar to Theorem \ref{thm:KapKleiner} for relatively hyperbolic pairs.

\begin{thm}[\cite{HaulRHBoundary}] \label{thm:Haul}
Let $(G, \P)$ be a rigid relatively hyperbolic pair with $\partial(G, \P)$ 1-dimensional and every $P \in \P$ one-ended. If each $P\in \P$ is tame, then $\partial(G, \P)$ is homeomorphic to one of the following:
\begin{enumerate}
    \item $S^1$
    \item $\mathcal{S}$, a Sierpinski carpet (planar)
    \item a Menger curve (non-planar)
\end{enumerate}
\end{thm}

By work of Dasgupta and Hruska \cite{DasguptaHruska}, the tameness condition on the peripherals can be dropped.

The theorem is based on the characterizations of global cut points in the Bowditch boundary. Just as Theorem \ref{thm:BowHyp} allows one to see splittings of hyperbolic groups via cut pairs in the boundary, further work of Bowditch shows that cut points in relatively hyperbolic boundaries correspond to splittings of the group over a subgroup of a peripheral group.

\begin{thm}[\cite{Bowperiph}]
Suppose $(G, \P)$ is a 1-ended relatively hyperbolic pair. If $G$ admits a peripheral splitting, then $\partial(G, \P)$ contains a global cut point.
\end{thm}

\begin{thm}[\cite{HaulRHBoundary}]
Suppose $(G, \P)$ is a 1-ended relatively hyperbolic pair with tame peripherals. If $\partial(G, \P)$ has a global cut point, then there exists a peripheral splitting of $(G, \P)$.
\end{thm}

Just as with Theorem \ref{thm:Haul}, tameness of peripherals can be dropped \cite{DasguptaHruska}.

Thus, a 1-ended relatively hyperbolic pair admits a peripheral splitting if and only if the Bowditch boundary contains a global cut point (if and only if there is a parabolic fixed point). See Figure \ref{fig:tree-of-circles} for an example of global cut points and a peripheral splitting.
Note that the Gromov boundaries of hyperbolic groups do not have cut points.

Cut points in Bowditch boundaries can cause exotic phenomena which prevent $G$ from being a Kleinian group. There are examples of relatively hyperbolic group pairs $(G, \mathcal{P})$ where the planar Bowditch boundary $\partial(G, \P)$ has cut points but where no peripheral structure on $G$ is virtually Kleinian or even virtually a manifold group \cite{HW1}.

This leads to the following question:

\begin{question}\label{question: vgfk}
When is a relatively hyperbolic group (virtually) a geometrically finite Kleinian group?
\end{question}

The Bowditch boundary of such a group must be planar, but having no cut points is not necessary. There are examples of geometrically finite Kleinian groups whose relatively hyperbolic boundary has cut points, such as a surface subgroup with accidental parabolics. Example \ref{example:TreeOfCircles} can be realized as a geometrically finite Kleinian group.

We have the following conjecture on the sufficient conditions \cite{HW1}.

\begin{conjecture}\label{conjecture: vgfk}
Let $(G,\P)$ be a relatively hyperbolic group pair. If $\partial(G,\P)$ is planar and has no cut points, then $G$ is virtually a geometrically finite Kleinian group.
\end{conjecture}

\begingroup
\sloppy
\hbadness=3000
\printbibliography
\endgroup

\bigskip

{\small
\noindent Michael Ben-Zvi, Bowdoin College \quad \texttt{mbenzvi@bowdoin.edu}

\noindent Jiayi Lou, Tufts University \quad \texttt{jiayi.lou@tufts.edu}

\noindent Genevieve S. Walsh, Tufts University \quad \texttt{genevieve.walsh@tufts.edu}
}

\end{document}